\documentclass{amsart}
\usepackage{amssymb}
\usepackage{amsmath}

\usepackage{a4}
\begin{document}
\newtheorem{theorem}{Theorem}[section]
\newtheorem{lemma}[theorem]{Lemma}
\newtheorem{remark}[theorem]{Remark}
\newtheorem{definition}[theorem]{Definition}
\newtheorem{corollary}[theorem]{Corollary}
\newtheorem{example}[theorem]{Example}
\newtheorem{assumption}[theorem]{Assumption}
\newtheorem{Rem}[theorem]{Remark} \newtheorem{sats}[theorem]{Theorem}
\newtheorem{prop}[theorem]{Proposition}
\newtheorem{lem}[theorem]{Lemma} \newtheorem{kor}[theorem]{Corollary}
\newcommand{\banm}{\begin{anm}}
\newcommand{\eanm}{\end{anm}}
\newcommand{\rank}{{\rm Im}}
\newcommand{\ind}{{\rm ind}}
\newcommand{\ord}{{\rm ord}}
\newcommand{\comp}{{\rm comp}}
\newcommand{\coker}{{\rm coker}}
\newcommand{\diag}{{\rm diag}}
\newcommand{\diam}{{\rm diam}}
\newcommand{\col}{{\rm col}}
\newcommand{\loc}{{\rm loc}}
\newcommand{\grad}{{\rm grad}}
\newcommand{\Dom}{{\rm Dom}}
\newcommand{\const}{{\rm const}}
\newcommand{\mes}{{\rm mes}}
\def\beq{\begin{eqnarray}}
\def\eeq{\end{eqnarray}}
\newcommand{\dr}{\frac{\partial}{\partial r}}
\newcommand{\nn}{\nonumber}
\newcommand{\nats}{\mbox{${\rm I\!N }$}}
\newcommand{\cal}{\mathcal}
\newcommand{\R}{\mathbb{R}}
\newcommand{\N}{\mathbb{N}}
\newcommand{\E}{\mathbb{E}}
\def\qedbox{\hbox{$\rlap{$\sqcap$}\sqcup$}}
\newcommand{\W}{{\mathaccent"7017 W}} \newcommand{\V}{{\mathaccent"7017 V}} \newcommand{\dist}{{\rm dist}} \font\pbglie=eufm10 \def\CD{C_{\mathcal{D}}} \def\CR{C_{\mathcal{R}}} \def\CDR{C_{\mathcal{D,R}}} \def\DD{\text{\pbglie D}} \makeatletter
  \renewcommand{\theequation}{%
   \thesection.\alph{equation}}
  \@addtoreset{equation}{section}
 \makeatother
\def\BB{\mathcal{B}}
\title[Wiener sausage and heat kernel asymptotics]
{Expected volume of intersection of Wiener sausages and
heat kernel norms on compact Riemannian manifolds with boundary}
\author{M. van den Berg and P. Gilkey}
\begin{address}{MvdB: Department of Mathematics, University of
Bristol, University Walk, Bristol,\newline\phantom{...a}BS8 1TW, U.K.}\end{address}
\begin{email}{M.vandenBerg@bris.ac.uk}\end{email}
\begin{address}{PG: Mathematics Department, University of Oregon, Eugene, OR 97403, USA}\end{address}
\begin{email}{gilkey@darkwing.uoregon.edu}\end{email}
\begin{abstract} Estimates are obtained for the expected volume of
intersection of independent Wiener sausages in Euclidean
space in the small time limit. The asymptotic behaviour of the weighted diagonal heat kernel norm on compact
Riemannian manifolds with smooth boundary is obtained in the small time limit.
\end{abstract}
\maketitle
\section{Introduction}

Let $B_1$, ..., $B_k$ be $k$ independent Brownian bridges in
$\mathbb{R}^m$ associated with the parabolic operator
$-\Delta+\frac{\partial}{\partial t},$ and with $B_i(0)=B_i(t)=0$
for $1\le i\le k$. Let $K$ be a compact, non-polar set in
$\mathbb{R}^m$. For $1\le i\le k$, let
$$
W^i_K(t)=\{B_i(s)+y:0\le s\le t,y\in K\}
$$
be the corresponding pinned Wiener sausages. These random sets are
Borel measurable with probability one, and we let
$$Z_{k,m}(t):=\mathbb{E}^1\otimes...\otimes\mathbb{E}^k\left[|\textstyle\prod_{i=1}^k\cap W^i_K(t)|\right]\,.$$
It is well known that $Z_{k,m}(t)$ is related to the virial
coefficients of a quantum system of obstacles $K$ at inverse
temperature $t$. For example, G. E. Uhlenbeck \cite{UB36}
calculated the asymptotic behaviour of $Z_{k,m}(t)$ as
$t\rightarrow\infty$ in the special case where $k=1$, $m=3$, and
$K$ is a ball. I. McGillivray \cite{McG00} obtained the full
asymptotic series as $t\rightarrow\infty$ for $k=1$, $m\ge 3$, and
$K$ an arbitrary, non-polar compact set. The special case where
$k=1$, $m=2$, and $K$ a ball was studied for $t\rightarrow\infty$
by M. van den Berg and E. Bolthausen \cite{BeBo97}.

In this paper, we analyse the asymptotic behaviour of $Z_{k,m}(t)$ as $t\rightarrow0$. Let $p_{\mathbb{R}^m}(\cdot,\cdot;\cdot)$
denote the heat kernel for $\mathbb{R}^m$ given by
$$
p_{\mathbb{R}^m}(x,y;t)=(4\pi t)^{-m/2}e^{-|x-y|^2/(4t)}\,,
$$
and let $p_{\mathbb{R}^m-K}$ denote the Dirichlet heat kernel for the open set $\mathbb{R}^m-K$. By the Feynman-Kac formula, we
have that for $x\in\mathbb{R}^m-K$,
\begin{equation}\label{eqn-1.a}
\mathbb{P}_x\left[B(s)\cap K=\emptyset,0\le s\le t\right]=(4\pi t)^{m/2}p_{\mathbb{R}^m-K}(x,x;t)\,,
\end{equation}
where $B(s)$,  $0\le s\le t$ is a Brownian bridge with
$B(0)=B(t)=x$. It is easily seen that
\begin{equation}\label{eqn-1.b}
0\le p_{\mathbb{R}^m-K}(x,x;t)\le p_{\mathbb{R}^m}(x,x;t)\,,
\end{equation}
and that
\begin{equation}\label{eqn-1.c}
Z_{k,m}(t)=(4\pi
t)^{km/2}\int_{\mathbb{R}^m}\left\{p_{\mathbb{R}^m}(x,x;t)-p_{\mathbb{R}^m-K}(x,x;t)\right\}^kdx\,,
\end{equation}
where the Dirichlet heat kernel is extended to all of $\mathbb{R}^m$ by putting
$$
  p_{\mathbb{R}^m-K}(x,y;t):=0\quad\text{for}\quad
  x\in K\quad\text{and or}\quad y\in K\,.
$$
This extends the equality of formula (\ref{eqn-1.a}) to all of $\mathbb{R}^m$.

It follows from the results in \cite{BrGi} that if $k=1$, $m\ge2$,
and $K$ is compact with smooth boundary, then there exists an
asymptotic series as $t\to 0$:
$$Z_{1,m}(t)=\sum_{j=0}^Jc_{1,j}t^{j/2}+O(t^{(J+1)/2})\,,
$$
where the coefficients $c_{1,j}$ have been computed for $j=0,...,4$. For example,
\begin{eqnarray*}
&&c_{1,0}=\int_K1dx,\quad
  c_{1,1}=\frac{\sqrt\pi}2\int_{\partial K}1dy,\quad\text{and}\quad
  c_{1,2}=\frac16\int_{\partial K}L_{aa}dy\,,
\end{eqnarray*}
where $dy$ is the surface measure on $\partial K$, $L_{aa}$ is the
trace of the second fundamental form on $\partial K$, and where
$\partial K$ is oriented with an outward orientation, i.e. the
inward normal is chosen on $\partial(M-K)$.

Let $\xi(t)$ be defined for $t>0$. We say that
$\xi(t)\sim\sum_{n\ge0}\xi_nt^{n/2}$ as $t\rightarrow0$ if
$$\xi(t)=\sum_{n=0}^J\xi_nt^{n/2}+O(t^{(J+1)/2})\quad\text{for any}\quad J\in\mathbb{N}\,.$$

The main results of this paper are the following two theorems:
\begin{theorem}\label{thm-1.1}
Let $K$ be a compact set in $\mathbb{R}^m$ ($m\ge2$) with smooth
boundary $\partial K$. Let $k\in\mathbb{N}$  be arbitrary. Then one has that
$Z_{k,m}(t)<\infty$ for all $t>0$, and
$Z_{k,m}(t)\sim\sum_{j\ge0}c_{k,j}t^{j/2}$ as $t\rightarrow0$
where,
\begin{eqnarray*}
&&c_{k,0}=\int_K1dx,\quad
  c_{k,1}=\frac12\left(\frac\pi k\right)^{1/2}\int_{\partial
  K}1dy,\\&&c_{1,2}=
 -\frac 13\int_{\partial K}L_{aa}dy,
\end{eqnarray*}
and for $k\ge2$
\begin{align*}
\begin{split}
    c_{k,2}=\biggl\{ &-\frac 1{2k}+\frac k2 \biggl(\frac\pi {4(k-1)^{3/2}}\\
    &-\frac 1{2k(k-1)}-\frac {\arctan (k-1)^{-1/2}}{2(k-1)^{3/2}} \biggr)\biggr\} \int_{\partial D} L_{aa} dy.
\end{split}
\end{align*}
\end{theorem}

\begin{theorem}\label{thm-1.2}
Let $M$ be a compact $m$-dimensional Riemannian manifold with smooth boundary $\partial M$. Let
$p_M(\cdot,\cdot;\cdot)$ denote the heat kernel for the Laplace-Beltrami operator acting on $L^2(M)$ with Dirichlet boundary
conditions on $\partial M$. Let $f$ be smooth on $M$ and let $k\in\mathbb{N}$.
\begin{enumerate}
\item For $t\rightarrow0$, $\displaystyle\int_M(p_M(x,x;t))^kf(x)dx\sim(4\pi
t)^{-km/2}\sum_{j=0}^\infty a_{k,j}t^{j/2}$.
\item $a_{k,0}=\displaystyle\int_Mfdx$.
\item
$a_{k,1}=\displaystyle\frac12\sqrt{\pi}\displaystyle\sum_{\ell=1}^k(-1)^\ell
\ell^{-1/2}\left(\begin{array}{l}k\\\ell\end{array}\right)\int_{\partial
M}fdy$.
\item $a_{k,2}=\displaystyle\frac
k6\int_ M f\tau dx$\newline
$\displaystyle-\frac 12
\displaystyle\sum^k_{\ell=1}\ell^{-1}\int_{\partial M}f^{(1)}
dy+\biggl\{-\displaystyle\frac k 6+\frac
12\sum^k_{\ell=1}\left(\begin{array}{ll}
    k\\
    \ell
   
\end{array}\right)(-1)^{\ell-1}$
\newline$\displaystyle\times\biggl\{(k-\ell)\left(\displaystyle\frac\pi{4\ell^{3/2}}-\displaystyle\frac
1{2\ell(\ell+1)}-\displaystyle\frac {\arctan\ell^{-1/2}}{2\ell^{3/2}}\right)+\displaystyle\frac 1\ell
    \biggr\}\biggr\}\displaystyle\int_{\partial M}fL_{aa}dy$

\end{enumerate}\noindent where $\tau$ is the
scalar curvature, and where $f^{(1)}$ is the normal derivative of
$f$ with respect to the inward unit normal vectorfield on
$\partial M$.\end{theorem}

It is possible to obtain the formulae for the $c_{k,j}, j\ge 1$ in
Theorem \ref{thm-1.1} directly from the formulae for the $a_{k,j},
j\ge 1$ in Theorem \ref{thm-1.2} by noting that the formulae in
Theorem \ref{thm-1.2} also hold for non-compact $m$-dimensional
Riemannian manifolds provided $\partial M$ is compact and smooth
and $f$ is smooth and has compact support. Specializing to the
case where $M=\mathbb{R}^m-K$, $f$ compactly supported and
identically equal to 1 in a neighbourhood of $\partial K$ and
noting that $ p_{\mathbb{R}^m}(x, x; t)=(4\pi t)^{-m/2}$ yields a
closed set of equations relating the $c_{k,j}$ and $a_{k,j}$ by
equating powers of $t$:
\begin{equation}\label{eqn-1.d}\begin{array}{l}
c_{k,j}=\displaystyle\sum_{\ell=1}^k(-1)^\ell\left(\begin{array}{l}k\\\ell\vphantom{{}_A}\end{array}\right)a_{\ell,j},\quad
j\in\mathbb{N}{}_{\vphantom{\vrule height 20pt}},\\
a_{k,j}=\displaystyle\sum^k_{\ell=1}(-1)^\ell \left(\begin{array}{c}
k\\
\ell
\end{array}\right)c_{\ell,j}, \quad j\in \mathbb{N}.\end{array}
\end{equation}

This paper is organized as follows. In Section \ref{sect-2}, we consider a half space and a suitable
localizing function $f$. The formulae for $a_{k,0}$ and
$a_{k,1}$ in Theorem \ref{thm-1.2} follow from Lemma \ref{lem-2.1}. One of the two boundary terms in
$a_{k,2}$ also follows from Lemma \ref{lem-2.1}. 

In Section
\ref{sect-3}, we study a planar region  in Euclidean space.  In Lemma \ref{lem-3.3}, we determine the
coefficient of $\int_{\partial D}L_{aa}dy$ in $c_{k,2}$. This completes the proof of Theorem \ref{thm-1.1}.  The
second boundary term in $a_{k,2}$ is determined in Lemma
\ref{lem-3.4} by computing the coefficient of $\int_{\partial
M}fL_{aa}dy$. This completes the computation of the boundary term in
$a_{k,2}$. It would be a tedious combinatorial exercise to show that these coefficients are consistent with
formula
\eqref{eqn-1.d}.

In Section \ref{sect-4}, we complete the proof of Theorem
\ref{thm-1.2} by determining  the interior term in $a_{k,2}$. This uses
the results of McKean and Singer \cite{MS67}. The arguments used in this section together with
existing formulae
\cite{Gi04} would suffice to determine $a_{k,\ell}$ for all $k$ and for
$\ell\le6$  if the boundary of $M$ is empty. The crucial difficulty in this paper, however, is the
determination of the boundary terms. These are not accessible by the methods used in Section \ref{sect-4}.

In Section
\ref{sect-5} we obtain results for the expected volume  of
intersection of Wiener sausages which are not pinned, and which
complement Theorem \ref{thm-1.1}.

\section{Half space calculations}\label{sect-2} Let $p_{\mathbb{R}^+}(x,x;t)$ be the Dirichlet heat kernel for the positive half
line $\mathbb{R}^+$.

\begin{lemma}\label{lem-2.1} Let $k\in\mathbb{N}$, and let $f$ be smooth. For $t\rightarrow0$, we have
\begin{eqnarray*}
&&\int_{\mathbb{R}^+}\left\{p_{\mathbb{R}^+}(x,x;t)\right\}^kf(x)dx\\
&\sim&(4\pi t)^{-k/2}\left\{
  \int_{\mathbb{R}^+}f(x)dx+\sum_{j=1}^\infty t^{j/2}f^{(j-1)}(0)\alpha_{k,j}\right\}
\end{eqnarray*}
where
$$\alpha_{k,j}=\frac12\sum_{\ell=1}^k(-1)^\ell\left(\begin{array}{l}k\\\ell\end{array}\right)
\Gamma\left(\frac j2\right)\Gamma(j)^{-1}\ell^{-j/2}\,.
$$
\end{lemma}

\begin{proof} We have that
\begin{eqnarray*}
&&p(x,x;t)=(4\pi t)^{-1/2}(1-e^{-x^2/t}),\\
&&f(x)\sim\sum_{j=0}^\infty\frac1{j!}f^{(j)}(0)x^j\quad\text{as}\quad x\rightarrow 0,\\
&&(1-e^{-x^2/t})^k=1+\sum_{\ell=1}^k(-1)^{\ell}\left(\begin{array}{l}k\\\ell\end{array}\right)e^{-\ell x^2/t},\\
&&\int_0^\infty x^je^{-\ell x^2/t}dx=\frac12\Gamma\left(\frac{j+1}2\right)\left(\frac
t\ell\right)^{(j+1)/2},\quad j>-1\,.
\end{eqnarray*}
Lemma \ref{lem-2.1} follows from these identities.\end{proof} Note that by
identity (1.45) of \cite{HWG72}
$$
\alpha_{k,2}=-\frac 12 \sum_{\ell=1}^k\ell^{-1}.
$$

\section{Computations for planar regions and proof of Theorem \ref{thm-1.1}}\label{sect-3}
\begin{lemma}\label{lem-3.1}
Let $K$ be a compact set in $\mathbb{R}^m$, and let $Z_{k,m}(t)$
be given by \eqref{eqn-1.c}. Then $Z_{k,m}(t)<\infty$ for all
$t>0$.
\end{lemma}
\noindent\emph{Proof.} Since $p_{{\mathbb{R}}^m}(x,y;t)\le (4\pi
t)^{-m/2}$ we have by \eqref{eqn-1.b}, \eqref{eqn-1.c} that

\begin{align}\label{eqn-3.a}
\begin{split}
Z_{k,m}(t)&\le (4\pi t)^{m/2}\int_{\mathbb{R}^m}\{p_{\mathbb{R}^m}(x,x;t)-p_{\mathbb{R}^m-K}(x,x;t)\}dx\\
&=(4\pi t)^{m/2} \int_{\R^m}\int_{\R^m}\{(p_{\R^m}(x,y;t/2))^2-(p_{\R^m-K}(x,y;t/2))^2\}dxdy\\
&\le 2^{(2+m)/2} \int_{\R^m}\int_{\R^m}\{p_{\R^m}(x,y;t/2)-p_{\R^m-K}(x,y;t/2)\}dxdy\\
&=2^{(2+m)/2}\int_{\R^m-K}\left\{1-\int_{\R^m-K}
p_{\R^m-K}(x,y;t/2)dy\right\}dx\\
&\qquad+2^{(2+m)/2}|K|,
\end{split}
\end{align} where we have used the semigroup property of the heat kernels for
$\R^m$ and $\R^m-K$ respectively. The first term in the right hand
side of \eqref{eqn-3.a} is the amount of heat in $\R^m-K$ at time
$t/2$ if $K$ is kept at fixed temperature $2^{(2+m)/2}$, and if
$\R^m-K$ has initial temperature 0. This term is finite by the
results of F.\ Spitzer \cite{FS12}. \qed

\medskip

A key ingredient in the proofs of Theorems \ref{thm-1.1} and
\ref{thm-1.2} is an estimate for the Dirichlet heat kernel
obtained by R.\ Lang  \cite{RL11} and  H.\ R.\ Lerche and D.
Siegmund \cite{HLDS10}. For related results see also  \cite{L68}.
Before we state their results (Theorem \ref{thm-3.2}) we introduce
some further notation. For an open, bounded and connected set $D$
in $\R^2$ we define the distance function $\delta:D\to (0,\infty)$
by $\delta(x)=\min\{|x-z|:z\in \partial D\}$. Suppose $\partial D$
is of class $C^3$. Then there exists $\varepsilon_D>0$ such that
for all $x\in D$ with $\delta(x)<\varepsilon_D$ there is a unique
$s\in \partial D$, depending on $x$, with $|x-s|=\delta(x)$. Hence
for all such $x\in D$ there is a smooth parametrization $x\to
(s,\delta)$ where $s\in \partial D$ is parametrized  by arc
length.
\begin {theorem}\label{thm-3.2}
Let $D$ be an open, bounded and connected set in $\R^2$ with $C^3$
boundary. Let $x\in D$ be such that $\delta(x)\le \varepsilon_D
/2$. Then as $t\to 0$

\begin{eqnarray}
p_D(x,x;t)&=&(4\pi t)^{-1}
\left\{1-e^{-\delta(x)^2/t}-L_{aa}(s)\delta(x)^2t^{-1/2}\int^\infty_{\delta(x)t^{-1/2}}e^{-\eta^2}
d \eta\right\}\nonumber
\\&&\qquad+O(1),\label{eqn-3.b}
\end{eqnarray}
where the remainder $O(1)$ is uniform on $\{x \in D: \delta(x)\le
\varepsilon_D /2\}$,  $L_{aa}$ is the curvature of $\partial D$ at
$s$, and  where $\partial D$ is oriented by an inward unit vector
field.
\end{theorem}

\begin{lemma}\label{lem-3.3}
Let $D$ be as in Theorem \ref{thm-3.2}, and let $k\in \N$. Then as
$t\to 0$
\begin{equation}\label{eqn-3.c}
    (4\pi t)^k \int_D \{p_{\R^2} (x,x; t)-p_D(x,x;t)\}^k dx=\sum^2_{j=0} c_{k,j} t^{j/2}+o(t),
\end{equation}
where $c_{k,0}=|D|, c_{k,1}=\frac 12 (\frac \pi
k)^{1/2}\int_{\partial D} 1 ds$,
$$
    c_{1,2}=-\frac 13 \int_{\partial D} L_{aa} (s) ds,
$$
and for $k\ge 2$ \begin{align*}
\begin{split}
    c_{k,2}=\biggl\{ &-\frac 1{2k}+\frac k2 \biggl(\frac\pi {4(k-1)^{3/2}}\\
    &-\frac 1{2k(k-1)}-\frac{\arctan (k-1)^{-1/2}}{2(k-1)^{3/2}} \biggr)\biggr\} \int_{\partial D} L_{aa}(s) ds.
\end{split}
\end{align*}
\end{lemma}
\noindent\emph{Proof.}  Let $\varepsilon\in (0,\frac 12)$ and
denote $D_\varepsilon=\{x\in D: \delta(x)\le~ t^\varepsilon\}$.
Since
\begin{equation}\label{eqn-3.d}
p_{\R^2}(x,x;t)\ge p_D (x,x;t)\ge p_{\R^2}(x,x;t)-(\pi t)^{-1}
e^{-\delta(x)^2/(2t)}
\end{equation}
(see Theorem 1 and estimate (ii) in \cite{MvdB13}) we that for
$\varepsilon< \frac 12$ and $t\to 0$
\begin{equation}\label{eqn-3.e}
\int_{D\backslash D_\varepsilon} \{p_{\R^2}
(x,x;t)-p_D(x,x;t)\}^k dx =
O(e^{-t^{2\varepsilon-1}/3}).\end{equation}
It remains to
calculate the contribution from the set $D_\varepsilon$ to the
integral in the left hand side of \eqref{eqn-3.c}. Since the
$O(1)$ remainder estimate is uniform in $x$ we have that this
remainder estimate contributes at most $O(t^{1+\varepsilon})$ to
the integral in \eqref{eqn-3.c}. Hence for $t$ sufficiently small
we have  by \eqref{eqn-3.e}
\begin{align*}\begin{split} Z_{k,2} (t)&= \int_{\partial D} ds
\int^{t^\varepsilon}_0 dr \biggl\{e^{-r^2/t}\\&+ L_{aa}(s) r^2
t^{-1/2}\int^\infty_{rt^{-1/2}}e^{-\eta^2}d\eta\biggr\}^k(1-rL_{aa}(s))+O(t^{1+\varepsilon}).\end{split}\end{align*}
Since $\int^{t^\varepsilon}_0 dr r^4t^{-1}=O(t^{5 \varepsilon-1})$
we have that for $\varepsilon\in \left(\frac 25, \frac 12\right)$
the contribution from such a term is $o(t)$. Hence for $t\to 0$
\begin{eqnarray}
Z_{k,2}(t)&=&\int_{\partial D} ds \int^{t^\varepsilon}_0 dr
\biggl\{e^{-kr^2/t}+ke^{-(k-1)r^2/t}
L_{aa}(s)r^2t^{-1/2}\int^\infty_{rt^{-1/2}} e ^{-\eta^2}
d\eta\biggr\}\nonumber\\
&&\times(1-rL_{aa}(s))+O(t^{5\varepsilon-1}).
\label{eqn-3.f}
\end{eqnarray}
Furthermore for $\varepsilon< 1/2$ and $t\to 0$
\begin{eqnarray}\label{eqn-3.g}
&&\int_{\partial D} ds \int_0^{t^\varepsilon} dr e^{-kr^2/t} \sim\frac
    12 \left (\frac{\pi t}k \right)^{1/2} \int_{\partial D} ds,
\\&&
    -\int_{\partial D} ds \int^{t^\varepsilon}_0 dr e^{-kr^2/t}
    rL_{aa} (s) \sim-\frac t{2k}\int_{\partial D} ds L_{aa}(s),
\nonumber\\&&
    \int_{\partial D}ds L_{aa} (s) \int_0^{t^\varepsilon}
    dr r^2 t^{-1/2} \int^\infty_{rt^{-1/2}} e^{-\eta^2} d \eta \sim
    \frac t6 \int_{\partial D} ds L_{aa}(s),\nonumber
\end{eqnarray}
and for $k\ge 2$

\begin{align}\label{eqn-3.h}
\begin{split}
&\int_{\partial D} ds L_{aa}(s) \int^{t^\varepsilon}_0 dr k
e^{-(k-1)r^2/t}
r^2 t^{-1/2} \int^\infty_{rt^{-1/2}} e^{-\eta^2} d \eta\\
&\sim \frac {kt}2 \int_{\partial D}ds L_{aa} (s) \int^\infty_0 dr k e^{-(k-1)r^2/t} r^2 t^{-1/2} \int^\infty_{rt^{-1/2}} e^{-\eta ^2} d\eta\\
&= \frac {kt}2 \int_{\partial D} ds L_{aa}(s)\int^\infty_1 d\eta (\eta^2+k-1)^{-2}\\
&= \frac{kt}2\left\{\frac \pi{4(k-1)^{3/2}}-\frac 1{2k(k-1)}-\frac
{\arctan (k-1)^{-1/2}}{2(k-1)^{3/2}}\right\}\int_{\partial D}ds
L_{aa} (s).
\end{split}
\end{align}
It is easily seen that the term with $L_{aa}(s)^2$ in
\eqref{eqn-3.f} contributes $O(t^{3/2})$. Collecting the
contributions from (\ref{eqn-3.g} -- \ref{eqn-3.h}) we arrive at
the conclusion of Lemma \ref{lem-3.3}. \qed

\medskip

 \noindent\emph{Proof of Theorem
\ref{thm-1.1}.} Since the contribution to the asymptotic series
for $Z_{k,m}(t)$ comes from a $t^\varepsilon, \varepsilon < 1/2$
neighbourhood of $\partial K$ we can identify the coefficients
$c_{k,0}, c_{k.1}$ for $m=2$ by Lemma \ref{lem-3.3}. By standard
results of invariance theory they also hold for $m>2$ \cite{Gi04}.
Moreover $Z_{k,m}(t)$ has a full asymptotic series. The proof of
Theorem \ref{thm-1.1} is complete by Lemma \ref{lem-3.1}.\qed

\begin{lemma}\label{lem-3.4}
The coefficient of the term $\int_{\partial D} fL_{aa}dy$ in the
expression for $a_{k,2}$ in Theorem \ref{thm-1.2} is given by
$$
    -\frac k6+\frac 12\sum^k_{\ell=1}\left(\begin{array}{ll}
    k\\
    \ell
    \end{array}\right)(-1)^{\ell-1}\biggl\{(k-\ell)\left(\frac\pi{4\ell^{3/2}}-\frac 1{2\ell(\ell+1)}
-\frac{\arctan\ell^{-1/2}}{2\ell^{3/2}}\right)+\frac 1\ell
    \biggr\}.
$$
\end{lemma}
\noindent\emph{Proof.} It suffices to consider a planar region
with $f\equiv1$ as in the proof of Lemma \ref{lem-3.3}. Let
$\varepsilon\in (2/5, 1/2)$. Then by \eqref{eqn-3.d}

$$\int_{D-D_\varepsilon}dx\{p_D(x,x;t)\}^k \sim(4\pi t)^{-km/2} |D-D_\varepsilon|.$$
The $O(1)$ remainder in \eqref{eqn-3.b} contributes at most
$O(t^{-k+1+\varepsilon})$ to the integral over the set
$D_\varepsilon$. Hence for $2/5<\varepsilon<1/2$
\begin{align*}
\begin{split}
\int_{D_\varepsilon}& dx \{p_D(x,x;t)\}^k=(4\pi
t)^{-k}\int_{\partial D}ds\int_0^{t^\varepsilon}
dr (1-rL_{aa}(s))\\
&\times \biggl\{1-e^{-r^2/t}-r^2t^{-1/2}L_{aa}(s)\int^\infty_{rt^{-1/2}}e^{-\eta^2}d \eta\biggr\}^k+O(t^{-k+1+\varepsilon})\\
&= (4\pi t)^{-k}\int_{\partial D}ds \int_0^{t^\varepsilon}dr
(1-rL_{aa}(s))\biggl\{(1-e^{-r^2/t})^k\\&-k(1-e^{-r^2/t})^{k-1}r^2t^{-1/2}L_{aa}(s)\int^\infty_{rt^{-1/2}}e^{-\eta^2}
d\eta\biggr\}+O(t^{5\varepsilon-k-1}),
\end{split}
\end{align*}
by estimates similar to the ones in the proof of Lemma
\ref{lem-3.3}. For $2/5<\varepsilon<1/2$ and $t\to 0$ we have
\begin{eqnarray*}
&&(4\pi t)^{-k} \int_{\partial D}ds \int^{t^\varepsilon}_0 dr
(1-rL_{aa}(s))=(4\pi t)^{-k} |D_\varepsilon|,\\
&&(4\pi t)^{-k} \int_{\partial D} ds \int_0^{t^\varepsilon} dr \{(1-e^{-r^2/t})^k-1\}\\
&\sim&2^{-1}(4\pi
t)^{-k}(\pi t)^{1/2}\sum^k_{\ell=1}\left(\begin{array}{ll} k\\\ell
\end{array}\right)(-1)^\ell \ell^{-1/2},
\end{eqnarray*}
and
\begin{eqnarray}\label{eqn-3.i}
&&-(4\pi t)^{-k} \int_{\partial D} ds L_{aa}(s)\int_0^{t^\varepsilon} dr r \{(1-e^{-r^2/t})^k-1\}\\
&\sim&2^{-1}(4\pi t)^{-k} t\sum^k_{\ell=1}\left(\begin{array}{l}
k\\\ell
\end{array}\right)(-1)^{\ell-1}\ell^{-1}\int_{\partial
D}ds L_{aa}(s).\nonumber
\end{eqnarray}
Furthermore
\begin{eqnarray*}
    &&-k(4\pi t)^{-k} \int_{\partial D} ds L_{aa}(s)\int^{t^\varepsilon}_0 dr
    r^2t^{-1/2}\int^\infty_{rt^{-1/2}} e^{-\eta^2} d\eta\\& \sim&
    -\frac k6(4\pi t)^{-k} t\int_{\partial D} ds L_{aa}(s),\nonumber
\end{eqnarray*}
and so it remains to compute  ( for $k\ge 2$)
\begin{eqnarray}\label{eqn-3.j}
&&-k(4\pi t)^{-k} \int_{\partial D}
dsL_{aa}(s)\int_0^{t^\varepsilon} dr\\&&\qquad\times
\sum^{k-1}_{\ell=1}\left(\begin{array}{c} k-1\\\ell
\end{array}\right)(-1)^{\ell}e^{-\ell r^2/t}r^2t^{-1/2}\int^\infty_{rt^{-1/2}} e^{-\eta^2} d\eta\nonumber\\
&&\sim -k(4\pi t)^{-k} \int_{\partial D} dsL_{aa}(s)\int_0^\infty
dr\nonumber\\&&\qquad\times \sum^{k-1}_{\ell=1}\left(\begin{array}{c} k-1\\\ell
\end{array}\right)(-1)^{\ell}e^{-\ell r^2/t}r^2t^{-1/2}\int^\infty_{rt^{-1/2}} e^{-\eta^2} d\eta\nonumber\\
&&= -\frac k2(4\pi t)^{-k} t\int_{\partial D} ds
L_{aa}(s)\sum^{k-1}_{\ell=1}\left(\begin{array}{c} k-1\\\ell
\end{array}\right)(-1)^{\ell} \int^\infty_1(\eta^2+\ell)^{-2}d\eta\nonumber\\
&&=-\frac k2(4\pi t)^{-k} t\int_{\partial D}ds
L_{aa}(s)\sum^{k-1}_{\ell=1}\left(\begin{array}{c} k-1\\\ell
\end{array}\right)\nonumber\\
&&\qquad\times(-1)^{\ell}\biggl\{\frac \pi{4\ell^{3/2}}-\frac
1{2\ell(\ell+1)}-\frac{\arctan\ell^{-1/2}}{2\ell^{3/2}}\biggr\}.\nonumber
\end{eqnarray}
Collecting all of the above we have shown that the formulae for
$a_{k,0}, a_{k,1}$ and $a_{k,2}$ in Theorem \ref{thm-1.2} hold for
planar regions with $C^3$ boundary. In particular we have shown
that, by collecting the terms (\ref{eqn-3.i}-\ref{eqn-3.j}), the
coefficient of the term $\int_{\partial D} fL_{aa} dy$ in the
expression for $a_{k,2}$ is given by Lemma \ref{lem-3.4}. The scalar
curvature term in the expression for $a_{k,2}$ in Theorem
\ref{thm-1.2} (4) will follow from the results in Section
\ref{sect-4} below.\qed

\bigskip

\section{The interior terms}\label{sect-4}
In order to obtain the contributions from the interior of $M$ to the coefficients $a_{k,j}$ in Theorem
\ref{thm-1.2}, it suffices to consider Riemannian manifolds without boundary. The general setting is as follows.
Let
$D$ be an operator of Laplace type on a smooth vector bundle $V$ over a Riemannian manifold $M$ without boundary.
One can express
$D$ in terms of geometrical data in an invariant fashion as follows. There is a unique connection
$\nabla$ on $V$ and there is a unique endomorphism
$E$ of $V$ so that one can express $D$ in the Bochner formalism by writing:
$$D=-\{g^{ij}\nabla_i\nabla_j+E\}\,.$$
Let $p_M(x,x;t)$ be the heat kernel, and let $F$ be a smooth
endomorphism of
$V$. Theorem
\ref{thm-1.2} generalizes to this setting to yield the existence of a complete asymptotic expansion
$$\int_M\operatorname{Tr}_{V_x}\{(p_M(x,x;t))^kF(x)\}dx\sim(4\pi t)^{-km/2}\sum_{j=0}^\infty
a_{k,2j}(F,D)t^{j},$$ where the coefficients $a_{k,2j}$ are locally computable. There are no half integer
powers of $t$ since the boundary of $M$ is empty. Consequently $a_{k,2j+1}=0$.

If $k=1$, then we have a complete
asymptotic expansion
\begin{equation}\label{eqn-4.a}
p_M(x,x;t)\sim(4\pi t)^{-m/2}\sum_{\nu=0}^\infty e_{2\nu}(D)t^\nu\end{equation}
where the $e_{2\nu}$ are locally computable endomorphisms. Let $\tau:=R_{ijji}$ be the scalar curvature.  
Then one has the
following well known result (see, for example, the discussion in \cite{Gi04}):
\begin{lemma}\label{lem-4.1} Let $D$ be a operator of Laplace type on a closed Riemannian manifold.
\begin{enumerate}
\item $a_{1,n}(F,D)=\int_M\operatorname{Tr}(Fe_n(x,D))dx$.
\item $e_0(x,D)=\operatorname{id}$.
\smallbreak\item $e_2(x,D)=\frac16\{\tau\operatorname{id}+6E\}$.
\end{enumerate}
\end{lemma}

\begin{remark}We note that similar formulae are available for $e_4$ and $e_6$ \cite{Gi04}.\end{remark}

Raising the asymptotic expansion in formula (\ref{eqn-4.a}) to the $k^{\operatorname{th}}$ power and applying
Lemma
\ref{lem-4.1} then yields immediately the following:

\begin{corollary}\label{cor-4.2} Let $M$ be a closed Riemannian manifold.\begin{enumerate}
\smallbreak\item $a_{k,0}(F,D)=a_{1,0}(F,D)$.
\smallbreak\item $a_{k,2}(F,D)=ka_{1,2}(F,D)$.
\smallbreak\item $a_{k,4}(F,D)=ka_{1,4}(F,D)+\int_M\operatorname{Tr}\{F\frac{k(k-1)}{72}(\tau^2+12\tau E+36E^2)\}dx$.
\end{enumerate}
\end{corollary}

\section{Expected volume of intersection of independent Wiener sausages}\label{sect-5}
In this section we obtain results analogous to Theorem
\ref{thm-1.1}  for the unpinned Wiener sausage.  Let $\beta^1(s),
\dots \beta^k(s), \quad s\ge 0$ be independent Brownian motions in
$\R^m$, and let $K$ be a compact set in $\R^m$. The Wiener
sausages $S^1_K(t), \dots, S^k_K(t)$ are the random sets defined
by
$$
    S_K^i(t)=\{\beta^i(s)+y : \quad 0\le s\le t, y\in K\}, \quad i=1, \dots,
    k.
$$
These random sets are compact and  Borel measurable with
probability 1. Define
$$ Q_{k,m}(t)= \E^1_0\otimes\dots\otimes \E^k_0\left[|\bigcap^k_{i=1} S_K^i(t)|\right].$$
It is well known that
$$
    Q_{k,m}(t)=\int_{\R^m} dx \biggl\{1-\int_{\R^m-K}p_{\R^m-K}(x,y;t)
    dy\biggr\}^k,
$$
and that $Q_{k,m}(t)\le Q_{1,m}(t)<\infty$ (\cite{FS12}). We note that
$Q_{k,m}(t)$ has the following analytic interpretation. Let
$u:\R^m-K\times (0, \infty)\to \R $ be the unique weak solution of
$\Delta u=\frac {\partial u}{\partial t}$ with initial condition
$u(x;0)=0, x\in \R^m-K$ and boundary condition $u(x;t)=1, x\in
\partial K, t>0$. We extend $u$ to all of $\R^m\times(0, \infty)$, and
note that
$$
    u(x;t)=1-\int_{\R^m-K} p_{\R^m-K}(x,y;t) dy.
$$
Then
$$Q_{k,m}(t)=\Vert u(\cdot; t)\Vert^k_k.$$
The large $t$ behaviour of $Q_{k,m}(t)$ has been investigated for
$m=2$ in \cite{JFLG14} and for $m\ge 3$ in \cite{MvdB15}.

\begin{theorem}\label{thm-5.1}
Let $K$ be a compact set  in $\R^m$ with $C^\infty$ boundary. There
exists an asymptotic series as $t\to 0$
$$
   Q_{k,m}(t) \sim \sum^{\infty}_{j=0} b_{k,j} t^{j/2},
$$
where
\begin{equation}\label{eqn-5.a}
b_{k,0}=\int_K 1 dx,
\end{equation}
\begin{equation}\label{eqn-5.b}
    b_{k,1}=2^{2+k}\pi^{(2-k)/2}
    \Gamma((k+1)/2)\int^\infty_1 d\eta_1\dots\int^\infty_1d\eta_k(\eta^2_1+\dots +\eta_k^2)^{-(k+1)/2}\int_{\partial K}1 dy,
\end{equation}
\begin{align}\label{eqn-5.c}
\begin{split}
    b_{k,2}&=-2 ^{2+k}\pi^{(2-k)/2}(k-2) \Gamma((k+2)/2)\\&\times\int^\infty_1
    d\eta_1\dots \int^\infty_1 d\eta_k(\eta^2_1+\dots +\eta^2_k)^{-(k+2)/2}\int_{\partial
    K} L_{aa} dy.
    \end{split}
\end{align}
\end{theorem}
\noindent\emph{Proof}. It suffices to prove
(\ref{eqn-5.a}--\ref{eqn-5.c}) for a disk in $\R^3$. The general
case for compact $K$ with $C^\infty$ boundary follows by
invariance theory. For the disk $D^3$ centred at the origin it is
well known \cite{CaJa86} that
$$
    u(x;t)=2\pi^{-1/2}r^{-1}
    \int^\infty_{(r-1)/(2t^{1/2})}e^{-\eta^2}d\eta,\qquad r=|x|\ge 1.
$$
Hence for $K=D^3$ we have that
$$Q_{k,3}(t)=\frac{4\pi}3+4\pi\left(\frac 2{\pi^{1/2}}\right)^k\int_1^\infty dr r^{2-k}\left\{\int^\infty_{(r-1)/(2t^{1/2})}e^{-\eta^2} d\eta\right\}^k.$$
The main contribution to $Q_{k,3}(t), t\to 0$ comes from the
interval $[1, 1+t^\varepsilon]$, where $0<\varepsilon<1/2$.
Expanding $r^{2-k}$ near 1 yields
\begin{align*}
\begin{split}
Q_{k,3}(t)&= \frac{4\pi}3+ 4\pi \left(\frac 2{\pi^{1/2}}\right)^k\int^\infty_0 dr \{1+(2-k)r\}\left\{\int^\infty_{r/(2t^{1/2)}}e^{-\eta^2}d\eta\right\}^k+O(t^{3/2})\\
&=\frac{4\pi}3+8\pi\left(\frac 2{\pi^{1/2}}\right)^k\int^\infty_0 dr r^k\int^\infty_1d\eta_1\dots\int^\infty_1 d\eta_ke^{-r^2(\eta^2_1+\dots+ \eta^2_k)}t^{1/2}\\
&- 8\pi\left(\frac 2{\pi^{1/2}}\right)^k(k-2)\int^\infty_0 dr r^{k+1}\int^\infty_1d\eta_1\dots \int^\infty_1 d\eta_k e^{-r^2(\eta^2_1+\dots+\eta_k^2)}t\\
&+O(t^{3/2}),
\end{split}
\end{align*}
and the claim follows by Fubini's theorem.\qed

We note that the coefficient of $t$ vanishes for $k=2$. This jibes
with the fact \cite{MvdB16} that for any compact $K$ in $\R^m$,
$m\ge1$, and all $t>0$

$$
    Q_{2,m}(t)=2Q_{1,m}(t)-Q_{1,m}(2t).
$$
\section*{Acknowledgments}
The research of P. Gilkey was partially supported by the Max
Planck Institute for the Mathematical Sciences (Leipzig, Germany).
The research of M. van den Berg was supported by the London
Mathematical Society under Scheme 4, reference 4511.

\end{document}